\documentclass[12pt]{article}

\usepackage[colorlinks=true, pdfstartview=FitV, linkcolor=blue, citecolor=blue, urlcolor=blue]{hyperref}

\usepackage{amssymb,amsmath, amscd}
\usepackage{times, verbatim}
\usepackage{graphicx}
\usepackage[OT2,OT1]{fontenc}

\DeclareFontFamily{OT1}{rsfs}{}
\DeclareFontShape{OT1}{rsfs}{n}{it}{<-> rsfs10}{}
\DeclareMathAlphabet{\mathscr}{OT1}{rsfs}{n}{it}

\newtheorem{prop}[equation]{Proposition}

\DeclareMathOperator{\End}{End}

\DeclareMathOperator{\Hom}{Hom}
\DeclareMathOperator{\Lie}{Lie}
\DeclareMathOperator{\Gal}{Gal}

\DeclareMathOperator{\disc}{disc}

\DeclareMathOperator{\Res}{Res}

\def\R{{\mathbb R}}

\def\bC{{\mathbb C}}
\def\Q{{\mathbb Q}}
\def\Z{{\mathbb Z}}

\def\O{{\mathscr O}}
\def \mod {{\rm mod}}

\title{On the periods of abelian varieties}

\author{Benedict H. Gross}

\begin{document}
\maketitle

\tableofcontents

\section {Abstract}

In this expository paper, we review the formula of Chowla and Selberg for the periods of elliptic curves with complex multiplication, and discuss two methods of proof. One uses Kronecker's limit formula and the other uses the geometry of a family of abelian varieties. We discuss a generalization of this formula, which was proposed by Colmez, as well as some explicit Hodge cycles which appear in the geometric proof.

\section{Looking for a thesis}

In my third year of graduate school at Harvard I was still looking for a thesis topic. John Tate had suggested a  problem on $p$-adic Galois representations, but I couldn't see how to make any progress on it. Fortunately for me, Neal Koblitz and David Rohrlich had arrived at Cambridge as BP assistant professors, and I started to talk to them about their work. Rohrlich showed me how the periods of eigen-differentials on the Fermat curve $F(d)$ of exponent $d$
could be explicitly calculated, using the values of Euler's gamma function at rational arguments with denominator $d$.  At the time, I was reading Andr\'e Weil's book ``Elliptic Functions according to Eisenstein and Kronecker'' \cite{W}. Weil ended with a proof of the Chowla-Selberg formula, which yields an expression, using the values of the gamma function at rational arguments with denominator $d$, for the periods of an elliptic curve with complex multiplication by an order in the imaginary quadratic field of discriminant $-d$. 

The similarity of these two expressions led me to wonder if there was any connection between them. In this direction, I was able to identify some products of gamma values at rational arguments $a/d$, which were the periods of certain differential forms of degree $n = \phi(d)/2$ on a factor $J$ of the Jacobian of the Fermat curve $F(d)$, and could also be viewed (via the Chowla-Selberg formula) as the periods of forms of the same degree on the product abelian variety $A^n$, where $A$ was an elliptic curve with complex multiplication by the imaginary quadratic field of discriminant $-d$. If I could relate the periods of these forms on the two abelian varieties geometrically, it would give an independent proof of the implication of the Chowla-Selberg formula for elliptic periods. (Their original proof used techniques from analytic number theory, including Kronecker's limit formula.) At the time I didn't see how this identification of periods could be made, but thought it might be a good topic to investigate further.

After I found the relation between the periods of higher degree forms, I  came across a paper \cite{W2} that Weil had recently written, which indicated that he was thinking about similar questions (and was way ahead of me). Since I heard that he was coming to Harvard to speak at Lars Ahlfors’ $70^{th}$ birthday conference, I wrote Weil a note asking if we could get together. He suggested that I meet him at the University guest house, and we ended up taking a long walk through the back streets of Cambridge. Weil described his own years as a graduate student in Rome, where had spoken with Vito Volterra about period integrals.

In fact, periods were the subject of Weil's talk at the Ahlfors conference the next day. He remained skeptical of William V.D. Hodge's famous conjecture describing the cohomology classes of algebraic cycles, and suggested that it might be time to search for a counter-example. David Mumford had found some interesting candidates in the Hodge ring of abelian varieties with complex multiplication, and Weil observed that these Hodge classes actually existed on a continuous family of abelian varieties, whose members had endomorphisms by an order in an imaginary quadratic field \cite{W3}. 

About half-way into this talk, with the single-minded focus that only a graduate student could possess, I realized that the techniques Weil was using could be modified slightly to give a geometric proof of the Chowla-Selberg formula for elliptic integrals. The continuous family of abelian varieties he highlighted (over a base which would now be called a Shimura variety associated to a group of unitary similitudes) had one fiber isomorphic to the factor $J$ of the Fermat Jacobian and another fiber isomorphic to the product $A^n$ of elliptic curves with complex multiplication. I could guess that the integrals of a differential form of higher degree on this family were constant (the form was horizontal for the Gauss-Manin connection), and that would give a geometric proof of the period identity. 

I spent the next few weeks learning the algebraic geometry that was necessary to write this all up, and sent the first draft of my argument in a letter to Weil (which eventually became the paper \cite{G}). He responded with a letter of encouragement, containing some suggestions for further work. The summer after Weil's lecture at Harvard, I attended a conference on automorphic forms and $L$-functions held in Corvallis, Oregon. All of the experts in the field had gathered there, and once they realized this, quite sensibly decided to lecture to each other. We graduate students present were frequently lost, and turned to Jean-Pierre Serre for help. After one of our remedial sessions, Serre mentioned that he had heard from Weil about my work on periods. He asked if I would like to meet Pierre Deligne, who was also thinking about these questions. After a quick introduction, the three of us sat down to talk.

I began by saying that I had found a new proof of the period implication of the Chowla-Selberg formula, using some techniques from algebraic geometry. Deligne immediately asked, in all seriousness, if I had proved the Hodge conjecture. I replied that I would be delighted to hear that I had done so, as I was still looking for a thesis topic (and felt that a proof of the Hodge conjecture would probably be sufficient). He then asked me to explain what I had actually done. After about fifteen minutes I had gone through the argument above, and we all agreed that I hadn't proved the Hodge conjecture. But Deligne thought that my argument could be extended to yield something in that direction. A few weeks later, he sent me a handwritten note of three pages outlining his proof of a fundamental theorem: that all Hodge cycles on abelian varieties are absolutely Hodge \cite{D}. So that was what I had been doing!
(For an illuminating general discussion of the Hodge conjecture, see \cite{D2}.)

About the same time, Serre asked me a question about Hecke characters for imaginary quadratic fields, and referred me to his paper with Tate  \cite{ST}, which gives an elegant treatment of the algebraic Hecke characters associated to abelian varieties with complex multiplication. I started to talk with Tate about it, and one thing led to another. So I ended up writing my thesis on the arithmetic of elliptic curves with complex multiplication \cite{G2}, rather than on their periods. 

In this expository paper, I will try to pull the two topics together. I will begin by reviewing the original analytic proof of the Chowla-Selberg formula. I will then introduce the elliptic curves $A = A(p)$ in my thesis, with complex multiplication by the integers of $k =\Q(\sqrt{-p})$. These curves are defined over the Hilbert class field $H$ of $k$, and are isogenous over $H$ to all of their Galois conjugates \cite{G2}.  I will show how the Chowla-Selberg formula gives information on the periods of $A$ at the complex places of $H$, and will reinterpret that result in terms of the Faltings height of $A$ \cite{F}. I will then describe a beautiful conjecture of Pierre Colmez \cite{C},
which gives a formula for the Faltings height of a general abelian variety with complex multiplication by the ring of integers of a CM field in terms of  logarithmic derivatives of Artin $L$-functions at $s=0$, and will summarize the recent progress that has been made in that direction.

I will end by explaining Deligne's construction of a motive of rank $2$ and weight $n$ with complex multiplication by an imaginary quadratic field $k$ from an abelian variety of dimension $n$ with endomorphisms by $k$ (which is an abstraction of the geometric argument in my Chowla-Selberg paper). I will work this construction out for the abelian variety $B = B(p) = \Res_{H/k} A(p)$ of dimension $h(-p)$ as well as for a factor $C$ of the Jacobian of the Fermat curve of exponent $p$, which has dimension $(p-1)/2$. The comparison of these two motives (they differ only by a Tate twist) yields a Hodge class in the middle cohomology of the product variety $B\times C$. In the simplest non-trivial case, when $p = 7$, $B$ is the elliptic curve $A(7) = X_0(49)$ with affine equation
$$y^2 + xy = x^3 - x^2 - 2x - 1$$
and $C$ is the Jacobian of the hyperelliptic curve of genus $3$ with affine equation
$$z^2 - z = t^7.$$
In this case, the abelian variety $B \times C$ has dimension $4$ and has a Hodge class of type $(2,2)$. Is there a codimension $2$ cycle on $B \times C$ which is defined over $\Q$ and has this class in cohomology?

\section{The Chowla-Selberg formula}

Let $k$ be an imaginary quadratic field with discriminant $-d$, ring of integers $\O$, class number $h$, and unit group of order $w$. Let $\mathfrak a_i$ be ideals of $\O$ which represent the distinct ideal classes. We fix an embedding of $k$ into $\bC$, so the ideals $\mathfrak a_i$ give lattices and the quotients $\bC/\mathfrak a_i$ correspond to the $h$ distinct complex elliptic curves with complex multiplication by $\O$.  Let $\Delta$ be the function on lattices in $\bC$ corresponding the the usual cusp form of weight $12$, and let $\Gamma(x)$ be Euler's Gamma function. Then the Chowla- Selberg formula is the equality \cite{CS}
$$\prod \Delta(\mathfrak a_i ) \Delta(\mathfrak a_i^{-1}) = (2\pi/d)^{12h} \prod\Gamma(a/d)^{6w\epsilon(a)}$$
where the first product is taken over the distinct ideal classes (the product $\Delta(\mathfrak a) \Delta(\mathfrak a^{-1})$ depends only on the ideal class of $\mathfrak a$)  and the second product is taken over the elements $0 < a < d$ which are relatively prime to $d$. Finally
$$\epsilon: (\Z/d\Z)^* \rightarrow \{\pm1\}$$
is the quadratic character which describes the splitting of primes in $k$, by quadratic reciprocity. 

The analytic proof involves the computation of the logarithmic derivative of the zeta function of $k$ at the point $s = 0$ in two different ways. Recall that the zeta function of $k$ is the Dirichlet series given by the sum over all non-zero ideals $\mathfrak a$ of the ring $\O$ of integers of $k$
$$\zeta_k(s) = \sum (N\mathfrak a)^{-s}$$
where the norm $N \mathfrak a$ of an ideal is equal to its index in the ring $\O$. This series converges in the half-plane where the real part of $s$ is greater than one. It has a meromorphic continuation to the entire complex plane, with a simple pole at $s = 1$ and no other singularities.

First, we may write the zeta function of $k$ as a sum of $h$ partial zeta funtions $\zeta_k(\mathfrak a_i,s)$, which are defined by taking the partial sum over the ideals $\mathfrak a$ of $\O$ in the same class as $\mathfrak a_i$. Kronecker's limit formula \cite[pg 73]{W} gives the first two terms in the Taylor expansion of $\zeta_k(\mathfrak a_i,s)$ at the point $s= 0$
$$\zeta_k(\mathfrak a_i,s) = -\frac{1}{w} - \frac{1}{12w}\log(\Delta(\mathfrak a_i)\Delta(\mathfrak a_i^{-1}))s + O(s^2)$$
Hence
$$\zeta_k(s) = -\frac{h}{w} - \frac{1}{12w} \sum \log(\Delta(\mathfrak a_i)\Delta(\mathfrak a_i^{-1}))s +O(s^2)$$
and  
$$d\log\zeta_k(s)|_{s=0} = \frac{1}{12h}\sum \log(\Delta(\mathfrak a_i)\Delta(\mathfrak a_i^{-1})).$$
 
On the other hand, we may also write the zeta function of $k$ as the product of the Riemann zeta function $\zeta(s)= \ \sum n^{-s}$ and the Dirichlet $L$ function $L(\epsilon,s) = \sum \epsilon(n) n^{-s}$ associated the the character $\epsilon$:
$$\zeta_k(s) = \zeta(s) L(\epsilon,s).$$
This identity was obtained by Dirichlet when the real part of $s$ is greater than $1$, by identifying the terms in the Euler product, using quadratic reciprocity. It then holds for all $s$ by analytic continuation.

From the product, it follows that
$$d\log\zeta_k(s) = d\log \zeta(s) + d\log L(\epsilon,s).$$
These logarithmic derivatives at $s = 0$ can be calculated from Lerch's expansion of the Hurwitz zeta function (with $0 < x \leq 1$) \cite[pg 59-60]{W}:
$$H(x,s) = \sum_{n = 0}^{\infty}\frac{1}{(n+x)^s} =  (1/2 - x) + \log (\Gamma(x)/\sqrt{2\pi})s + O(s^2).$$
Taking $x = 1$, we find that $\zeta(0) = -1/2$ and $d\log \zeta(s)|_{s = 0} = \log(2\pi).$
On the other hand, summing over $0 < a < d$ with $a$ prime to $d$ we find
$$L(\epsilon, s) = d^{-s} \sum \epsilon(a)H(a/d,s).$$
Since $\sum \epsilon(a) = 0$, we find $L(\epsilon,0) = -\sum \epsilon(a) (a/d)$ and $\zeta_k(0) = 1/2\sum \epsilon(a)(a/d))$. Comparing the last formula with the formula for $\zeta_k(0)$ obtained by summing the partial zeta functions gives Dirichlet's famous class number formula
$$h = -(w/2) \sum_{0 < a < d} \epsilon(a) (a/d).$$

We also obtain the following formula for the logarithmic derivative 
$$d\log L(\epsilon,s)|_{s=0} = (w/2h)\sum_{0 < a < d} \epsilon (a) \log \Gamma(a/d) - \log d.$$
Adding this to the logarithmic derivative of $\zeta(s)$ at $s = 0$, we get a second expression for the logarithmic derivative of $\zeta_k(s)$:
$$d\log \zeta_k(s)|_{s=0} = \log (2\pi) - \log d + (w/2h) \sum_{0 < a < d} \epsilon(a) \log \Gamma(a/d).$$
Comparing this with the one obtained via Kronecker's limit formula, multiplying both sides by $12h$, and exponentiating gives the Chowla-Selberg formula:
$$\prod \Delta(\mathfrak a_i ) \Delta(\mathfrak a_i^{-1}) = (2\pi/d)^{12h} \prod\Gamma(a/d)^{6w\epsilon(a)}$$

\section{Elliptic periods}

In this section, I will describe what the Chowla-Selberg formula yields on the periods of elliptic curves with complex multiplication in the simplest case, when the discriminant $d$ of the imaginary quadratic field is a prime.  Let $p$ be a prime number with  $p \equiv 3~ (\rm mod~4)$. Let $k = \Q(\sqrt{-p})$ be the imaginary quadratic field of discriminant $-p$ and let $\O$ be the ring of integers of $k$. We will also assume that $p > 3$, so that the group of units $\O^* = \langle \pm 1 \rangle$ and $w = 2$. %For future reference, we note that the imaginary quadratic field $k$ embeds as a subfield of the $p^{th}$ cyclotomic field $K = \Q(\mu_p)$. Indeed, the Galois group of $K/\Q$ is canonically isomorphic to $(\Z/p\Z)^*$, via the action on $p^{th}$ roots of unity, and $k$ is the subfield fixed by the subgroup $(\Z/p\Z)^{*2}$ of squares. An explict embedding is given by the Gauss sum
%$$g = \sum_{(\Z/p\Z)^*} (a|p) \zeta^a = \sum_{(\Z/p\Z)^*} \chi(a) \zeta^a$$
%where $\chi(a) = (a|p) : (\Z/p\Z)^* \rightarrow \langle \pm 1 \rangle$ is the quadratic residue symbol and $\zeta$ is a $p^{th}$ root of $1$ in $K$. Indeed, Gauss proved that $g^2 = -p$, so $k = \Q(g) \subset K$.

In this simple case, the character $\epsilon(a) = (a|p)$ is just the quadratic residue symbol and the Chowla- Selberg formula states:
$$\prod_{i = 1}^h \Delta(\mathfrak a_i ) \Delta(\mathfrak a_i^{-1}) = (2\pi/p)^{12h} \prod_{a = 1}^{p-1}\Gamma(a/p)^{12\epsilon(a)}.$$

To interpret this as a result on elliptic periods, we will introduce the elliptic curves $A(p)$. Let $H$ be the Hilbert class field of the imaginary quadratic field $k$, which is an abelian extension of $k$ with Galois group isomorphic to the ideal class group of $\O$. Elliptic curves $A$ with complex multiplication by $\O$ over $H$ are determined up to isomorphism by two invariants \cite[\S9]{G2}: an algebraic Hecke character
$$\chi_A: I_H^* \rightarrow k^*$$
on the id\`eles of $H$ whose restriction to the principal id\`eles $H^*$ is given by the norm, and the modular invariant $j(A) = j(\tau) \in H$. Here $\tau$ is a point in the upper half-plane which is a root of an integral quadratic polynomial $ax^2 + bx + c$ with discriminant $b^2 -4ac = -p$. The character $\chi_A$ determines the isogeny class of $A$ over $H$ and the invariant $j(A)$ determines the isomorphism class of $A$ over $\overline{\Q}$. 

The elliptic curve $A = A(p)$ has invariant $j(A) = j((1 + \sqrt{-p})/2)$, which generates the subfield $F = H \cap \R$,  and character $\chi_A$ of conductor $(\sqrt{-p})$. On id\`eles $b = (b_v)$ where the local components at places $v$ dividing $p$ are units which are congruent to $1$, let $\mathfrak b$ be the corresponding fractional ideal of $H$, and let $\mathfrak a = N_{H/k}(\mathfrak b)$. Then $\mathfrak a$ is principal, and $\chi_A(b)$ is the unique generator of $\mathfrak a$ which is a square modulo $\sqrt{-p}$. Since the character $\chi_A$ is equivariant for the action of $\Gal(H/\Q)$ on both $I_H$ and $k^*$, the curve $A = A(p)$ is isogenous to all of its conjugates over $H$ \cite[\S11]{G2}. The curve $A$ descends to its field of moduli $F = \Q(j(A))$, where it defines an isogeny class containing two isomorphism classes. We specify an isomorphism class by insisting that $A$ has a minimal Weierstrass model over the ring of integers of $F$ with discriminant $\Delta= -p^3$ \cite{G3}. The N\'eron differential $\omega$ determined by this model is well-defined up to sign, and the Chowla-Selberg formula gives a formula for the product of its periods over the complex places of $H$:
$$\prod_{v | \infty}~~ \int_{A_v(H_v)} |\omega_v \wedge \overline{\omega_v}| = (2\pi/p)^h \prod_{a = 1}^{p-1} \Gamma(a/p)^{\epsilon(a)}.$$
Taking logarithms of both sides and dividing by $h$, we find the equivalent formula
$$\frac{1}{h} \log~ (\prod_{v | \infty}~~ \int_{A_v(H_v)} |\omega_v \wedge \overline{\omega_v}| )= \log 2\pi - \log p + \frac{1}{h} \sum \chi(a) \log \Gamma(a/p) =   d\log \zeta_k(s)|_{s=0}.$$

We now derive the formula for the product of periods from the Chowla-Selberg formula. Since $H = \Q(\sqrt{-p}, j(A))$, we have a unique embedding $v_1: H \rightarrow \bC$ which maps $\sqrt{-p}$ to a complex number with positive imaginary part, and $j(A)$ to the real number $j((1 + \sqrt{-p})/2)$. For each class $\mathfrak a$ let $\sigma_{\mathfrak a}$ be the corresponding element in the Galois group of $H/k$, using the Artin reciprocity law, and let $v_{\mathfrak a}: H \rightarrow \bC$ be the embedding defined by $v_{\mathfrak a} = v_1 \circ \sigma^{-1}_{\mathfrak a}$. Then the complex period lattice of $\omega_{v_{\mathfrak a}}$ has the form $\Omega_{\mathfrak a}\cdot \mathfrak a$, for a non-zero complex number $\Omega_{\mathfrak a}$. The period integral at the place $v = v_{\mathfrak a}$ is then given by 
$$\int_{A_v(H_v)} |\omega_v \wedge \overline{\omega_v}| = \Omega_{\mathfrak a}\cdot\overline{\Omega_{\mathfrak a}} \cdot N(\mathfrak a) \sqrt p.$$
On the other hand, since $\Delta$ is a modular form of weight $12$, we have
$\Delta(\Omega_{\mathfrak a} \cdot \mathfrak a) = \Omega_{\mathfrak a}^{-12} \Delta(\mathfrak a)$,
and the value of $\Delta$ on the period lattice of $\omega_{v_{\mathfrak a}}$ is given by
$\sigma^{-1}_{\mathfrak a}(\Delta(\omega)) = \sigma^{-1}_{\mathfrak a}(-p^3) = -p^3.$
This gives the formula
$-p^3 \cdot\Omega_{\mathfrak a}^{12}=\Delta(\mathfrak a).$
Since $\overline{\Delta(\mathfrak a)} = \Delta(\overline{\mathfrak a}) = N(\mathfrak a)^{-12} \Delta(\mathfrak a^{-1})$ we conclude that
$$\Omega_{\mathfrak a}^{12} \cdot \overline{\Omega_{\mathfrak a}}^{12} \cdot N(\mathfrak a)^{12} \cdot p^6= \Delta(\mathfrak a) \Delta(\mathfrak a^{-1}) $$

Now take the product over all embeddings and use the formula for
$\prod \Delta(\mathfrak a ) \Delta(\mathfrak a^{-1})$ given by Chowla and Selberg. We find that
$$\prod_{v | \infty}~~ \int_{A_v(H_v)} |\omega_v \wedge \overline{\omega_v}|^{12} = (2\pi/p)^{12h} \prod_{a = 1}^{p-1} \Gamma(a/p)^{12\epsilon(a)}.$$
Since both sides are positive real numbers, we may take the $12^{th}$ roots to obtain the stated formula on elliptic periods. For generalizations to elliptic curves with complex multiplication by non-maximal orders, see \cite{K}, \cite{NT}.

We will need another result on the periods of regular differentials $\omega$ on $A$ over $H$, when integrated over rational $1$-cycles on the curve at the completions $H_v$ \cite[Thm 21.2.2]{G2}. Define an equivalence relation $a \sim b$ on non-zero complex numbers if the ratio $a/b$ lies in $k^*$. Then we have
$$\prod_{v |\infty} ~~\int_{\gamma_v} \omega_v ~~\sim~~(2 \pi i)^{-m}\prod_{\epsilon(a) = +1} \Gamma(a/p)$$
where $\gamma_v$ is any non-trivial $1$-cycle in the rational homology of $A$ over the complex completion $H_v$ and
$$m = \sum_{\epsilon(a) = +1} a/p  =  \frac{p-1}{4} - \frac{h}{2}$$ 
The product of $\Gamma$ values and the above sum is taken over
the quadratic residues in $(\Z/p\Z)^*$, and $a$ is the unique representative of the class which lies between $1$ and $p$. %Since
%$$\Gamma(x+1) = x\Gamma(x)$$
%the choice of representative does not change the equivalence class of the value $\Gamma(a/p)$.

Let $B = B(p) = \Res_{F/\Q} A(p)$. Then $B$ is an abelian variety of dimension $h$ over $\Q$ which has complex multiplication over $k$ by a CM field $E$ which contains $k$ \cite[\S 15]{G2}. The field $E$ is generated over $k$ by certain $h^{th}$ roots of ideals which are $h^{th}$ powers. Let $\omega_B$ be a non-zero regular differential of degree $h$ on $B$ over $k$. Then $\omega_B$ is unique up to scaling by $k^*$ and it follows from the above that the non-zero periods are given up to equivalence by
$$\int_{\gamma_B} \omega_B ~ \sim ~ (2 \pi i)^{-m}\prod_{\epsilon(a) = +1} \Gamma(a/p)$$
where $\gamma_B$ is any $h$-cycle in the rational homology of $B$.

\section{Colmez's conjecture for the Faltings height}

The formula we obtained for the product of the periods of $(1,1)$ forms $\int_{A_v(H_v)} |\omega_v \wedge \overline{\omega_v}|$, using the values of the Gamma function at rational numbers with denominator $p$, can be used to compute the Faltings height of the abelian variety $A = A(p)$ over $\overline{\Q}$. To define this height, we let $\alpha = \sqrt{-p}$~ in $H$ and pass to the quadratic extension $K = H(\sqrt\alpha)$ where the elliptic curve $A$ has good reduction everywhere. Let $\mathscr A$ be the N\'eron model of $A$ over $K$ and let $\Omega(\mathscr A)$ be the projective $\O_K$-module of N\'eron differentials. In this special case, the projective module $\Omega(\mathscr A)$ is free, and generated by $\omega' = \sqrt\alpha \cdot \omega$. The Faltings height $h(A)$ of $A$ over $\overline{\Q}$ is then defined as
$$h(A) = \frac{-1}{[K:\Q]} \sum_w \log\int_{A_w(\bC)} |\omega'_w \wedge \overline{\omega'_w}|$$
where the sum is taken over the complex places $w$ of $K$. 

Since there are two complex places $w$ of $K$ above each complex place $v$ of $H$, and
$$\int_{A_w(\bC)} |\omega'_w \wedge \overline{\omega'_w}| = \sqrt p \int_{A_v(\bC)} |\omega_v \wedge \overline{\omega_v}|$$
we obtain the formula
$$h(A) = \frac{-1}{2h} \sum_v \log\int_{A_v(\bC)} |\omega_v \wedge \overline{\omega_v}| ~~-~~ \frac{1}{4}\log p$$
Combining this with the Chowla-Selberg formula for $A(p)$, we find that
$$h(A) = -\frac{1}{2} d\log \zeta_k(s)|_{s=0} -\frac{1}{4} \log p = -\frac{1}{2}d\log L(\epsilon,s)_{s = 0} - \frac{1}{4}\log p -\frac{1}{2}\log(2\pi) .$$

While looking for a product formula for periods, analogous to the classical
product formula for algebraic numbers, Pierre Colmez was led to a beautiful generalization of the above result. He conjectures a precise formula, expressing the Faltings height of an abelian variety with complex multiplication in terms of the logarithmic derivatives of Artin $L$-functions at $s = 0$ \cite{C}. The amazing idea of relating the periods of abelian varieties with complex multiplication by an abelian field to the logarithmic derivatives of Dirichlet $L$-functions at $s = 0$ is due to Anderson \cite{A}. Colmez was able to extend this work and establish his product formula in the abelian case (with an assist from Obus \cite{O} at the prime $2$).  Substantial progress on the general Colmez conjecture was recently made by Andreatta, Goren, Howard, and Madapusi-Pera \cite{AGHM}, and by Yuan and Zhang \cite{YZ}. These authors establish an average version of the conjectural formula, which I will describe below.

Let $E$ be a CM field of degree $2g$, with ring of integers $\O_E$ and totally real subfield $E^+$. Let $B$ be an abelian variety of dimension $g$ with complex multiplication by $\O_E$, defined over $\overline{\Q}$. Then $B$ is defined over a finite extension $K$ of $\Q$, which is contained in $\overline{\Q}$. We will assume that $K$ is large enough so that all endomorphisms of $B$ are defined over $K$, $K$ contains the normal closure of $E$, and the abelian variety $B$ has everywhere good reduction over $K$ (cf. \cite{ST}). The CM type $\Phi$ of $B$ is defined as the set of embeddings $E \rightarrow K$ which result from the diagonalization of the action of $\O_E$ on the $g$ dimensional tangent space $\Lie(B/K)$. If $c$ is complex conjugation on $E$, then the union $\Phi \cup \Phi\circ c$ is the set of all embeddings of $E$ into $K$. Hence there are $2^g$ possible CM types $\Phi$, for each CM field $E$ of degree $2g$. The automorphism group of the normal closure of $E$ acts on the finite set of CM types for $E$ by composition.

Let $\mathscr B$ be the N\'eron model of $B$ over the ring of integers $\O_K$ of $K$, and let $\Omega = \det(\Lie(\mathscr B)^{\vee})$ be the projective $\O_K$ module (of rank $1$) of N\'eron differentials. Choose a non-zero element $\omega \in \Omega$, and define the Faltings height $h(B)$ of $B$ by the fomula
$$[K: \Q] \cdot h(B) = -\sum_v \log\int_{B_v(\bC)} |\omega_v \wedge \overline{\omega_v}| + \log\#(\Omega/\O_K\omega)$$
where the sum is taken over the complex places $v$ of $K$. Faltings shows that the height is independent of the choice of differential $\omega$ and the field of definition $K$ of $B$, and Colmez shows that it depends only on the CM type $\Phi$ of $B$ (in fact, it depends only on the Galois orbit of the CM type), so we can denote it $h(E,\Phi)$. His conjecture gives a precise formula for the height $h(E, \Phi)$ in terms of logarithmic derivatives of Artin $L$-series at $s = 0$. For $g = 1$ this is what we obtained from the Chowla-Selberg formula, and Tonghai Yang \cite{Y} has resolved the case when $g = 2$. The average version of the Colmez conjecture, which was recently proved, is the simpler statement that
$$\frac{1}{2^g} \sum_{\Phi} h(E,\Phi) = -\frac{1}{2} d \log L(V,s)_{s = 0} - \frac{1}{4} \log f(V) - \frac{g}{2} \log 2\pi$$
where the sum is taken over all the possible CM types for $E$, $V$ is the $g$ dimensional representation of the Galois group of $\Q$ induced from the non-trivial quadratic character of $E/E^+$, $L(V,s)$ is its Artin $L$-function, and $f(V) = (-1)^g\disc E/\disc E^+$ is its conductor.

A generalization of Colmez's conjecture to the logarithmic derivatives of Artin $L$-functions at all negative integers was proposed by Maillot and Roessler \cite{MR}. There has also been recent progress on the conjecture that I made with Deligne in \cite{G}, giving the periods of varieties acted on by automorphisms of finite order in terms of rational values of the $\Gamma$ function \cite{Fr} \cite{MR2}.

\section{Deligne's  motive}

In this section, we recall Deligne's construction of a motive of rank $2$ and weight $n$ from an abelian variety $A$ of dimension $n$ with complex multiplication by an imaginary quadratic field $k$. This construction uses two facts. First, the higher de Rham, Betti, and \'etale cohomology groups of an abelian variety $A$ are given by the exterior powers of the first cohomology: $H^n(A) = \wedge^nH^1(A)$. Second, the exterior powers of a $k$ vector space $V$ embed in a natural way as direct factors of its exterior powers over $\Q$.

Let $A$ be an abelian variety of dimension $n$ over $\Q$, which has endomorphisms by the imaginary quadratic field $k$ over the extension field $k$. By this we mean that there is a homomorphism
$$k \rightarrow \End_k(A) \otimes \Q.$$
These endomorphisms act linearly on the $k$ vector space $\Lie(A/k)$, which decomposes as $u$ copies of the identity embedding of $k$ and $v$ copies of the conjugate embedding, with $u + v = n$.
Deligne constructs a subspace of dimension $2$ in the middle cohomology group $H^n(A)$ of $A$, in all cohomology theories (Betti, deRham, Hodge, and $\ell$-adic). This should be a motive $M = M(A)$ in the sense of Grothendieck; for us it will suffice that we can define its cohomological realizations. In that sense, the definition of $M$ is given by
$$M =  \wedge^n_kH^1(A) \subset \wedge^nH^1(A) = H^n(A).$$ 
For example, suppose $V = H_B^1(A)$ is Betti cohomology, which is a vector space of dimension $2n$ over $\Q$ with an action of $k$. Then the exterior power $\wedge^n_kV = M_B$ (which has dimension $1$ over $k$ and dimension $2$ over $\Q$) embeds as a direct factor of the $\Q$-vector space $\wedge^nV = H_B^n(A)$: it is the subspace of $\wedge^nV$ on which elements $\alpha$ in the group $k^*$ act by the two characters $\alpha^n$ and $\overline{\alpha}^n$.

The motive $M$ has
weight $n$ and rank $2$ over $\Q$ and has complex multiplication by $k$ over the extension field $k$. The Hodge numbers of $M_B \otimes \bC$ are $(u,v) + (v,u)$. Indeed, the two characters $\alpha^n$ and $\overline{\alpha}^n$of $k^*$ appear in bi-degrees $(u,v)$ and $(v,u)$ respectively. This shows that the Hodge decomposition is algebraic, and occurs over $k \subset \bC$. The periods of $M$ are the integrals of a deRham class $\omega_M$ of type $(u,v)$ over $k$ over rational classes in the dual of the one dimensional $k$-vector space $M_B$. These define an equivalence class in $\bC^*/k^*$, which determines the Hodge structure of $M$.

The $\ell$-adic realization $M_{\ell}$ of $M$ is the induced representation of the character
$$\rho_{M,\ell}: \Gal(\overline{k}/k) \rightarrow (k\otimes\Q_{\ell})^*$$
which comes from the determinant of the Galois represention on the $k \otimes \Q_{\ell}$-vector space
$H^1_{\ell}(A)$. All of these $\ell$-adic characters should come from a single algebraic Hecke character on the id\`eles $I_k$ of $k$:
$$\psi_M: I_k \rightarrow k^*$$
which is equivariant for complex conjugation and has algebraic part $z \rightarrow z^u\overline{z}^v$. The process of passing from an algebraic Hecke character to a compatible family of $\ell$-adic characters is described below.

We will now work out the realizations of the motive $M$ for the abelian variety $B = B(p)  = \Res_{F/\Q}A(p)$ of dimension $n = h(-p)$ over $\Q$. In this case, the action on $\Lie (B/k)$ is by $n$ copies of the standard embedding of $k$, so the Hodge numbers of $M_B \otimes \bC$ are
$(n,0) + (0,n)$. The periods of $M_B$ are given by the integrals of the non-zero regular differential $\omega_B$ defined over $k$ over classes $\gamma$ in the rational homology:
$$\int_{\gamma_B} \omega_B ~ \sim ~ (2 \pi i)^{-m}\prod_{\chi(a) = +1} \Gamma(a/p),$$
where $m = \sum_{\epsilon(a) = +1} \langle a/p \rangle = (p-1)/4 - h/2$

We can determine the $\ell$-adic character $\rho_{M,\ell}$ as follows. The abelian variety $B$ has complex multiplication over $k$ by the CM field $E$ of degree $2h$ \cite[\S15] {G2}. This determines an algebraic Hecke character $\psi_{M(B)}: I_k^* \rightarrow E^*$ whose algebraic part is the standard embedding of $k^*$ into $E^*$ \cite[\S8]{G2}. The representation of $\Gal(\overline{\Q}/k)$ on $H^1_{\ell}(B)$ can be described as follows \cite{ST}. Let $f_{\ell}$ be the embedding $(k \otimes \Q_{\ell})^* \rightarrow (E \otimes \Q_{\ell})^*$ and let 
$$\psi_{\ell}: I_k \rightarrow (E \otimes \Q_{\ell})^*$$
be defined by $\psi_{\ell}(a) = \psi_{M(B)}(a) \cdot f_{\ell}(a_{\ell})^{-1}$. Then $\psi_{\ell}$ is trivial on $k^*$. Since $\psi_{\ell}$ is continuous and its image is totally disconnected, it is also trivial on the connected component of the id\`ele class group, and the group of connected components of the id\`ele class group is isomorpic to the abelianized Galois group of $k$ (via the inverse of the Artin reciprocity law). This gives a homomorphism $\psi_{\ell}: \Gal(\overline{k}/k) \rightarrow (E \otimes \Q_{\ell})^*$, and the representation of this Galois group on $H^1_{\ell}(B)$ decomposes as the direct sum of the $2h$ characters obtained by the different embeddings of $(E \otimes \Q_{\ell})$ into $\overline{\Q_{\ell}}$. Its determinant  is therefore given by the composition of $\psi_{\ell}$ with the norm $N$ from $(E \otimes \Q_{\ell})^*$ to
$(k \otimes \Q_{\ell})^*$. 

It follows that the $\ell$-adic character $\rho_{M,\ell}$ corresponds to the algebraic Hecke character
$$\psi_M = N \circ \psi_B: I_k \rightarrow k^*$$
which is equivariant for complex conjugation.
The Hecke character $\psi_M = N \circ \psi$ has algebraic part given by the map $\alpha \rightarrow \alpha^h$ and conductor $(\sqrt{-p})$. For an ideal $\frak a$ which is prime to $p$, the value $\psi_M(\frak a) = N \circ \psi_B(\frak a)$ is the unique generator $\beta$ of the ideal $(\frak a)^h$ in $k$ which is congruent to a square $(\mod \sqrt{-p})$. This follows from the fact that $\psi_B(\frak a)$ is the unique $h^{th}$ root of $\beta$ in the CM field $E$. 

The induced $2$-dimensional $\ell$-adic representations $M_{\ell}$ of $\Gal(\overline{\Q}/\Q)$ correspond to a holomorphic modular form of weight $h + 1$ for the group $\Gamma_0(p^2)$, which is a newform with integral Fourier coefficients.

\section{A factor of the Fermat Jacobian}

We will now compute Deligne's motive $M = M(J)$ for a factor $C$ of the Jacobian of the Fermat curve of exponent $p$. This factor is the Jacobian of the quotient curve $X = X(r,s,t)$ with affine equation ~\cite{GR}
$$y^p = x^r(1-x)^s,$$
where $(r,s,t)$ is a triple of integers with $0< r,s,t < p$ and $r+s+t = p$. Like the Fermat curve, the curve $X$ is an abelian cover of $\mathbb P^1$, which is ramified only at $\{0,1,\infty\}$. The map from the Fermat curve $u^p + v^p = 1$ to the curve $X$ is given by $(u,v) \rightarrow (x,y) = (u^p,u^rv^s)$.

For a rational number $x$ with $n \leq x < n+1$ we let $\langle x \rangle = x - n$, so $0 \leq \langle x \rangle < 1$. The genus of $X$ is $n = (p-1)/2$ and for $a \in (\Z/p\Z)^*$ the differentials
$$\omega_a = x^{\langle ar/p\rangle - 1}(1-x)^{\langle as/p \rangle - 1}dx$$
give a basis for the first de Rham cohomology over $\Q$. The curve $X$ has an automorphism of order $p$ over $K = \Q(\mu_p)$, given by $(x,y) \rightarrow (x,\zeta.y)$, where $\zeta$ is a primitive $p^{th}$ root of unity. The differential $\omega_a$ is an eigenvector with eigenvalue $\zeta^a$. Moreover, $\omega_a$ is of the first kind if and only if \cite[\S 2]{G} \cite[pg 815]{W}
$$\langle ar/p \rangle + \langle as/p \rangle + \langle at/p \rangle = 1.$$
The periods of $\omega_a$ over any $1$-cycle $\gamma$ in the rational homology of $X$ have the form \cite[Appendix]{G2} \cite[pg 815]{W}
$$\int_{\gamma} \omega_a = I(\gamma)^{\sigma_a}B(\langle ar/p\rangle, \langle as/p \rangle)$$
where $I(\gamma)$ is an element in $K = \Q(\mu_p)$ which depends on the rational cycle $\gamma$ and $\sigma_a$ is the automorphism of $K$ that maps $\zeta$ to $\zeta^a$. Finally, $B(u,v)$ is Euler's beta function, defined by the integral
$$B(u,v) = \int_0^1 x^{u-1} (1-x)^{v-1} dx.$$

The differentials $\omega_a$ correspond to eigenforms of degree one on the Jacobian $C = C(r,s,t)$ of the curve $X$. This abelian variety has dimension $n =(p-1)/2$; over the extension field $K = \Q(\mu_p)$ it has complex multiplication by the ring $\Z[\mu_p]$ of integers in $K$. The CM type of $C$ consists of the elements
$$\Phi = \{a \in (\Z/p\Z)^*: \langle ar/p \rangle+ \langle as/p \rangle + \langle at/p \rangle = 1\}.$$ We will henceforth assume that $p \equiv 3~ (\mod ~4)$ and that $ p > 3$, and will consider $C$ as an abelian variety of dimension $n$ over $\Q$ which has endomorphisms by the imaginary quadratic field $k = \Q(\sqrt{-p})$ over $k$. Indeed, $k$ is the unique quadratic subfield of $K = \Q(\mu_p)$ in this case. The subgroup of the Galois group $(\Z/p\Z)^* = \Gal(K/\Q)$ which fixes $k$ consists of the squares 
$$\Gal(K/k) = \{a \in (\Z/p\Z)^*: \epsilon(a) = +1\}.$$ 

Our objective is to compute the various realizations of Deligne's rank $2$ motive $M = M(C)$, and compare them to the motive $M(B)$ studied in the previous section.
The Hodge numbers of $M_B \otimes \bC$ are $(u,v) + (v,u)$, where
%Consider the form of degree $n = (p-1)/2$ on $J$ which is the wedge product (in some order) of the one forms $\omega_a$ over the quadratic residues $a$
%$$\omega_J = \Lambda_{\chi(a) = +1} \omega_a$$
%This has pure Hodge type $(u,v)$ with 
$$u = \#\{a \in \Phi: \epsilon(a) = +1\}$$
$$v = \#\{a \in \Phi: \epsilon(a) = -1\}$$
Indeed, the CM type $\Phi$ gives the action of $K$ on $\Lie(C/K)$, so the action of the subfield $k$ is by $u$ copies of the standard embedding and $v$ copies of the conjugate embedding.
Clearly, $u + v = n$. Using Dirichlet's class number formula, one can show \cite[Lemma 4]{G} that
$u - v = h.\epsilon(r,s,t) $ 
with $\epsilon(r,s,t) = \epsilon(r) + \epsilon(s) + \epsilon(t)$.

In \cite[Thm 2]{G} I computed the periods of  the form $\omega_C$ of type $(u,v)$ on $C$ from the periods of the $1$-forms $\omega_a$. All the cycles of degree $n = (p-1)/2$ come from products of $1$-cycles, and
$$\int_{\gamma_1.\gamma_2\ldots \gamma_n}\omega_C =  \prod_{\epsilon(a) = +1} B(\langle ar/p\rangle, \langle as/p \rangle) \det((I(\gamma_i)^{\sigma_a})).$$
The latter determinant lies in $k$, as applying an element of the Galois group of $K/k$ to the determinant permutes the columns. Since this Galois group is abelian of odd order, this permutation is even and the determinant is unchanged. Hence
$$\int_{\gamma} \omega_C~\sim~ \prod_{\epsilon(a) = +1} B(\langle ar/p\rangle, \langle as/p \rangle).$$
Using Euler's relation between the beta and gamma functions, and his functional equation for the gamma function:  
$$B(x,y) = \Gamma(x)\Gamma(y)/\Gamma(x+y) ~~~~~\Gamma(x)\Gamma(1-x) = \pi/\sin(\pi x)$$
we find after some calculation that this product of beta values is equivalent to
$$(2 \pi i)^{-n} \prod_{\epsilon(a) = +1}\Gamma(ar/p)\Gamma(as/p)\Gamma(at/p)$$
where we recall that $n = (p-1)/2$.

If $\epsilon(r) = +1$ we have 
$$\prod_{\epsilon(a) = +1} \Gamma(ar/p) ~\sim~\prod_{\epsilon(a) = +1} \Gamma(a/p)$$
whereas if $\epsilon(r) = -1$ we have 
$$\prod_{\epsilon(a) = +1} \Gamma(ar/p) ~\sim~ (2\pi i)^n~/~\prod_{\epsilon(a) = +1} \Gamma(a/p).$$
Since the same holds for $s$ and $t$, when $\epsilon(r,s,t) = +1$ we find the simple formula
$$\int_{\gamma} \omega_C ~\sim~\prod_{\epsilon(a) = +1} \Gamma(a/p).$$ 
 This gives the periods of $M(C)$ when $\epsilon(r,s,t) = +1$. In this case, $u - v = h$ and $u + v = n = (p-1)/2$.
Hence $v = (p-1)/4 - h/2$ and
$$\int_{\gamma} \omega_C \sim \int_{\gamma} \omega_B \cdot (2\pi i)^v$$
Hence the Betti, de Rham, and Hodge realizations of $M(C)$ are Tate twists of the corresponding realizations of $M(B)$.

%We can also consider the $d$-form on $J$ defined by
%$$\nu_J =\Lambda_{\chi(a) = -1} \omega_a$$
%This is pure of Hodge type $(n,m)$. When $\epsilon(r,s,t) = +1$ a similar argument shows that
%$$\int_{\gamma} \nu_J ~\sim~\prod_{\chi(a) = -1} \Gamma(a/p)$$ 

Next, we compute the $\ell$-adic realization of $M(C)$. The representation of $\Gal(\overline{\Q}/k)$ on $H^1_{\ell}(C/k)$ is induced from an abelian representation of $\Gal(\overline{\Q}/K)$ corresponding to the Hecke character $\psi_C= \psi(r,s,t)$ given by Jacobi sums \cite[\S1 \S3]{GR}. This algebraic Hecke character is an Galois equivariant homorphism on the id\`eles of $K$
$$\psi_C:I_K \rightarrow K^*$$
whose algebraic part is the map of tori $K^* \rightarrow K^*$ determined by the CM type $\Phi(r,s,t)$ of $C$.

The determinant of the induced representation is then given by the transfer of the inducing character, as the sign of the permutation representation of the odd abelian group $\Gal(K/k)$ is trivial. By class field theory, the transfer of the inducing character is given by the restriction of the algebraic Hecke character $\psi_C$ to the id\`eles of the subfield $k$ \cite[Ch XIII]{S2}. This restriction gives a Galois equivariant homorphism
$$\psi_{M(C)}: I_k \rightarrow k^*$$
whose algebraic part is the homomorphism $z \rightarrow z^u\overline{z}^v$. Since the conductor of the Jacobi sum Hecke character $\psi_C$ is either $(1-\zeta)$ or $(1-\zeta)^2$, its restriction $\psi_{M(C)}$ to $I_k$ has conductor dividing $(\sqrt{-p})$. Since $u + v = (p-1)/2$ is odd, the conductor is equal to $(\sqrt{-p})$. A bit more analysis shows that when $\epsilon (r,s,t) = +1$, we have
$$\psi_{M(C)} = \psi_{M(B)} \cdot N^v$$
in $\Hom(I_k,k^*)$, where $N$ is the id\`elic norm.
Since the $2$-dimensional representation of $\Gal(\overline{\Q}/\Q)$ on $H_{\ell}(M(C))$ is induced from this character of $\Gal(\overline{\Q}/K)$, we can also identify the $\ell$-adic realization of $M(C)$ with a Tate twist of the $\ell$-adic realization of $M(B)$. Summing up, we have shown

\begin{prop}
Assume that $\epsilon(r,s,t) = \epsilon(r) + \epsilon(s) + \epsilon(t) = +1$. Then $u - v = h$ and the Deligne motives $M(B)$ and $M(C)$ of rank $2$ associated to the abelian varieties $B = B(p)$ and $C = C(r,s,t)$ differ by a Tate twist
$$M(C)(m)  = M(B)$$
with 
$m =  \sum_{\epsilon(a) = +1} \langle a/p\rangle = (p-1)/4 - h/2$.
\end{prop}

A similar analysis yields the identity $M(C)(m) = M(B)$ when $\epsilon(r,s,t) = -1$.

\section{A Hodge class}

In this section, we assume that $\epsilon(r,s,t) = \pm1$. We recall that $C = C(r,s,t)$ is the Jacobian of the curve $X(r,s,t)$ with equation $y^p = x^r(1-x)^s$ and $B = \Res_{F/\Q}A(p)$.  Let $d = (p-1)/2 + h = \dim(B \times C)$. Note that $d$ is even, as both $(p-1)/2$ and $h$ are odd. 

The rank $4$ motive $M(B) \otimes M(C)(d/2)$, occurs as a submodule of $H^d(B \times C)(d/2)$ by the Kunneth decomposition. It has Hodge numbers $(h,-h) + ((0,0) + (0,0) + (-h,h)$. Since the rank $2$ motives $M(B)$ and $M(C)$ are both symplectic and differ by a Tate twist, it follows that the tensor product $(M(B) \otimes M(C)(d/2)$ contains the rank $2$ Artin motive $\Q + \Q(\epsilon)$. In particular, there is a Hodge class in the middle cohomology of $B \times C$ which is defined over $\Q$.

Are these Hodge classes algebraic, or do they give counter-examples to the Hodge conjecture? For example, when $p \equiv 7 ~(\mod~8)$ the triple $(r,s,t) = (1,1,p-2)$ has $\epsilon (r,s,t) = +1$. In this case, the curve $X = X(r,s,t)$ is hyperelliptic, with affine equation
$$z^2 - z = t^p$$
Is there an algebraic cycle of codimension $d/2 = (p-1)/4 + h/2$ which is defined over $\Q$ and gives this Hodge class in the middle cohomology of the abelian variety $B \times C$?

\def\noopsort#1{}
\providecommand{\bysame}{\leavevmode\hbox to3em{\hrulefill}\thinspace}

Benedict H. Gross\\
Leverett Professor of Mathematics, Emeritus\\
Harvard University\\

\medskip
Harvard Department of Mathematics\\
1 Oxford St.\\
Cambridge, MA 02138\\
USA\\
\medskip

gross@math.harvard.edu

\end{document}